\documentclass [11pt,reqno]{article}
\usepackage[leqno]{amsmath}
\usepackage{paralist, amsthm,amsfonts,amsmath,amssymb,graphics,hyperref,pstricks,graphicx,multicol}
\usepackage{amsmath}    
\usepackage{graphicx}   
\usepackage{verbatim}   
\usepackage{hyperref}   
\usepackage{graphics}
\usepackage{amscd,amssymb}

\divide
\oddsidemargin by 2 
\makeatletter
\newcommand{\address}[2][]{%
  \ifx\@add@ress\@undefined\gdef\@add@ress{\par\par\bigskip}\AtEndDocument{\@add@ress}\fi
  \g@addto@macro\@add@ress{\bigskip\noindent{\small\scshape%
      \ifx#1\empty\else{\bfseries Address of #1:}\ \fi#2}\par\par}}

\renewenvironment{abstract}{\small\quotation\noindent
  {\bfseries \abstractname}}{\endquotation \par}
\newcommand{\footnotetextplain}[1]{\begingroup\def\@thefnmark{}%
  \long\def\@makefntext##1{\parindent 0pt\noindent ##1}\@footnotetext{#1}
  \endgroup}

\makeatletter \xdef\qedbuit{\qed}

\newcommand{\TeoremaAmbFinalMarcat}[1]{%
  \expandafter\gdef\csname end#1\endcsname{\qedbuit\@endtheorem}}

               {\hfill\rule{2.5mm}{2.5mm} \vspace{\parskip} } 

\theoremstyle{definition}
 \TeoremaAmbFinalMarcat{defi}
 \TeoremaAmbFinalMarcat{ex}
 \TeoremaAmbFinalMarcat{rem}
 \TeoremaAmbFinalMarcat{conv}

\newcommand{\start}[2]{\begin{#1}\label{#2}}

\def\@enum@{\list{\csname label\@enumctr\endcsname}%
           {\usecounter{\@enumctr}\def\makelabel##1{\hss\llap{##1}}
           \itemsep=2pt\parsep=0pt\topsep=3pt plus 1pt minus 1 pt}}

\newenvironment{numlist}{\enumerate[(1)]}{\endenumerate}

\newcommand{\bnl}{\begin{numlist}}
\newcommand{\enl}{\end{numlist}}

\begin{document}

\title{A combinatorial algorithm for visualizing representatives with minimal self-intersection}
\author{Chris Arettines}
\maketitle

\begin{abstract}

Given an orientable surface with boundary and a free homotopy class, we present a purely combinatorial algorithm which produces a representative of that homotopy class with minimal self intersection.

\end{abstract}

\section{Introduction}

	Given an orientable surface with boundary described by a surface word and a reduced word representing a free homotopy class(see Section \ref{sec:Basics}), we construct two sequences which encode the intersection structure of a representative curve(see Section \ref{sec:bo}).  We then examine these combinatorial sequences to see if they indicate the presence of bigons in the associated representative (see Section \ref{sec:fb}). If there are no bigons, then by a theorem of Hass and Scott, Theorem \ref{hs}, our representative must have minimal self-intersections.  If the sequences indicate the presence of a bigon, we must examine the combinatorial data to see if the bigon is \textit{improper} or \textit{proper}.  Proper bigons are essentially those that can be removed by a homotopy to reduce the total number of intersections(see Section \ref{sec:Basics} for a rigorous definition).  In the presence of a proper bigon, manipulations of the combinatorial data mimic the appropriate homotopy to reduce the total number intersections in the representative.  The algorithm then iterates until there are no more proper bigons to remove, which by Theorem \ref{hs}, indicates that our representative is minimal.  The heart of this paper is Section \ref{sec:rb}, which contains the theorems giving combinatorial classifications for proper and improper bigons.  This allows us to completely avoid the use of geometry, making the algorithm very simple to implement on a computer. 

An earlier draft of this paper was the author's undergraduate thesis at Stony Brook, under the guidance of Moira Chas.  The author is extremely thankful to her for the many insightful discussions, and for suggesting this interesting subject matter.  The author would also like to thank Tony Phillips and Dennis Sullivan for their helpful comments.

This algorithm has been programmed by the author and can be used at \href{http://www.math.sunysb.edu/~moira/applets/chrisApplet.html}{this link}. 

\section{Basics}
\label{sec:Basics}

\begin{figure}[http]
\center{\includegraphics[width=1.0in]{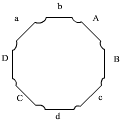}}
\caption{The surface word $abABcCdD$ and the associated polygon.  The unlabeled edges are curved to emphasize that they are part of the boundary.}
\label{sw}
\end{figure}

	We begin with a brief discussion of the basics and terminology used in this paper.  Given a set of symbols \textit{S} $=$ \{$s_1,s_2,\dots,s_{n},S_1,S_2,\dots,S_{n}$\}, we define a \textit{surface word} to be a cyclic sequence of these 2n symbols where each symbol in \textit{S} appears exactly once.  Given a surface word $w_1w_2...w_{2n}$ with each $w_i$ $\in$ \textit{S}, we associate to it a surface with boundary. Beginning with a polygon having 4n edges, we label the edges in a clockwise fashion as follows:  Choose one edge and label it $w_1$, leave the second edge unlabeled, label the third $w_2$, the fourth unlabeled, and so on.  We identify the edge labeled $s_i$ with the edge labeled $S_i$ in such a way as to preserve orientability.  This identification gives rise to an orientable surface with boundary.  In Figure \ref{sw}, we consider the set of symbols \{$a,b,c,d,A,B,C,D$\} with $a$ identified with $A$, $b$ identified with $B$, etc. (We will use this convention throughout the rest of this paper) The surface word $abABcCdD$ then represents the 2-torus with one boundary component.
	After being given a surface word in some alphabet \{$s_1,s_2,\dots,s_{n},S_1,S_2,\dots,S_{n}$\} , the fundamental group for the associated surface is then the free group generated by \{$s_1,s_2,\dots,s_{n}$\} and their inverses \{$S_1,S_2,\dots,S_{n}$\}.  Each free homotopy class is given by a unique reduced cyclic word in the fundamental group (Recall that a cyclic word is reduced if no element and it's inverse are adjacent for all cyclic permutations of the word).
	
	From a surface word $X$ and a reduced cyclic word $W$, we may construct a representative of the associated free homotopy class by first choosing distinct points along the appropriate edges of the polygon associated with $X$ that our curve will pass through and then proscribing which points should be connected by line segments(of course we must take care to ensure that we end up with a single closed curve).  Such a representative on a fundamental polygon is called a \textit{segmented representative}.  Aside from this pictorial description of a representative, we may also think of a representative as a map from the circle $S$ to our surface $M$.  This permits the following definition:

A map $f$ $S \rightarrow M$ is said to contain a $bigon$ if two conditions hold: (1) there exists two closed subarcs $A_1$ and $A_2$ of $S$ such that the endpoints of $A_1$ are mapped to the same points on M as the endpoints of $A_2$ are; (2) The images of the two arcs bound a topological disk when we choose appropriate lifts to the universal cover of $M$.  If the arcs $A_1$ and $A_2$ are disjoint, we call the bigon \textit{proper}, and call it \textit{improper} otherwise.

Not all representatives have the same number of self intersections.  We would like to have a systematic way of creating representatives with the minimal number of self intersections possible. This topic has been investigated in the past, and in fact there is the following result which seems intuitively obvious but is not trivial to prove:

\start{theo}{hs} (Hass and Scott, \cite{hassscott} ) If a representative $F$ does not have the minimal number of self intersections, then it must have a proper bigon.
\end{theo}

Once we have any starting representative, by the above theorem we simply need to look for proper bigons and homotope them away.  The algorithm to be described takes this "hands on" approach and translates it into a combinatorial method that a computer can easily perform.
\section{The Algorithm}
\label{sec:alg}
	\subsection{Basic Objects}
	\label{sec:bo}
	The input for the algorithm will be a surface word $X$ and a reduced cyclic word $W$.  Consider these fixed for the present discussion and suppose that a segmented representative in general position is given to us (i.e. all intersections occur in the interior of the fundamental polygon).  We now describe a way of associating to this representative, two combinatorial objects which encode the intersection structure.  Consider the first letter in $X$, say it is $a$, and look at the corresponding edge of the fundamental polygon.  Label the points along that edge through which the representative passes in a clockwise manner by $a_1,...,a_k$, where $k$ is the total number of appearances of $a$ and $A$ in $W$.  Since our curve is in general position, each point will be distinct.  We must also label the paired points on the $A$ edge in clockwise manner by $A_k,...,A_1$.  This is so that $a_i$ and $A_i$ really represent the same point on the surface.  Proceed clockwise around the edges of the polygon until all points where the curve intersects an edge are labeled.  Following the clockwise order around the polygon, we may sequence these points in a cyclic list, called a \textit{point list} (denoted P).  
	
\begin{figure}[http]
\center{\includegraphics[width=4.00in]{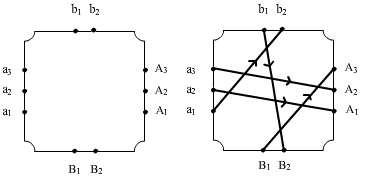}}
\caption{We start by placing the points in P, then connect them according to C to get a representative.}
\label{Alg1}
\end{figure}

	For example, if $X=abAB$ and $W=AAAbb$, then $P=a_1,a_2,a_3,b_1,b_2,A_3,A_2,A_1,B_2,B_1$ as in Figure \ref{Alg1}.  After choosing a labeling of points, we follow the orientation of the curve and sequence the pairs of points that are connected by line segments in another cyclic list, called a \textit{segment list} (denoted C). In Figure \ref{Alg1}, C $=$ $(B_1,A_3) (a_3,A_2) (a_2,A_1) (a_1,b_2) (B_2,b_1)$.  We call elements of P \textit{points} and elements of C \textit{word segments}, and we shall call identical lower and uppercase letters with the same index inverses of each other.  Hence $a_1^{-1}=A_1$, $B_1^{-1}=b_1$ and so on.  Note that while P will initially be the same for every representative of AAAbb, C will vary depending on how these points are connected by line segments.  The reader may also notice that the information in C is redundant - for example, a pair $(x,A_i)$ will always be followed by a pair whose first element is $a_i$.  However, the reader will see that this excess notation for C is convenient. 

	Conversely, just using P, C and X from above, we may construct a segmented representative.  Take the first element in P, say $a_1$, and choose a point on the $a$ edge of our polygon associated with X to label $a_1$.  This choice fixes our placement of the point labeled $A_1$ as well.  Proceed clockwise around the edges of the polygon, choosing and labeling points in the correct relative positions on the edges on and the correct sides.  Lastly, we use C to connect the appropriate points with line segments to complete the representative.  Since the precise geometric positioning of the points to connect with segments was arbitrary (only the relative positioning was important), our new representative may not be exactly identical to the original one we started with.  However, it will be shown that any two segmented representations constructed from the same P,C and X have the same number of intersections.

\subsection{Finding Bigons}
	\label{sec:fb}
	We may use P and C to determine if our representative has any bigons.  Let $C=W_0W_1...W_{(n-1)}$ where each $W_i=(w_{i_1},w_{i_2})$ is a word segment.  Consider two word segments $W_i$ and $W_j$.  The associated line segments in a representative intersect if and only if the pair of points $w_{i_1}$ and $w_{i_2}$ separate $w_{j_1}$ and $w_{j_2}$ in the cyclic list P (this corresponds to the line segments being transverse to each other). If this is the case, we say $W_i$ intersects $W_j$, and denote this combinatorial intersection $W_i$I$W_j$.  Immediately, we see that the number of intersections in a segmented representative determined by P and C depends only on P and C, and not on the specific geometric positioning of the points.  We summarize this discussion in the following proposition:
	
	\start{prop}{intersections} Intersections in a representative are in one-to-one correspondence with combinatorial intersections. \end{prop}

	Now, given an intersection $W_i$I$W_j$, we would like to determine if it is a vertex of a bigon.  To do this pictorially, we start at the intersection, and trace along pairs of line segments through repeated copies of the fundamental polygon in an attempt to identify the two "legs" of the bigon (see the example at the end of the paper to help clarify things). We try this for each of the four possible pairs of line segments until either the pair of line segments we were tracing intersect again, or the pair of segments don't intersect but split to different edges of the polygon .  In the first case, we have found a bigon, while in the second case, we know we may stop tracing in that direction, since in the universal cover, the lifts would lead to different fundamental regions.  In the case of surfaces without boundary, these lifts could still possibly meet up again to form a closed loop, but in the present case of surfaces with boundary, there is no possibility of the lifts meeting again.  If the tracings split in all four directions with no intersection encountered, we conclude that our initial point of intersection cannot possibly be part of a bigon.    
	
\begin{figure}[http]
\center{\includegraphics[width=2.0in]{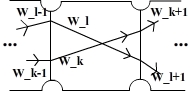}}
\caption{There are 4 directions to check in, but here, two of them split.}
\label{CombBigon}
\end{figure}

	This bigon finding procedure just described can easily be done combinatorially.  Starting with the intersection $W_k$I$W_l$, there are four possible pairs of word segments we may consider next. The four possible pairs of word segments we should look at are $W_{k+1}$ and $W_{l+1}$, $W_{k+1}$ and $W_{l-1}$,$W_{k-1}$ and $W_{l+1}$, and finally $W_{k-1}$ and $W_{l-1}$ (all indices taken mod $n$). See Figure \ref{CombBigon}.  To determine if the pair we are looking at intersects, we check the positions of the $w_{i_j}$s in $P$ as described above.  To determine if the pair splits to different sides, we also simply need to check the labelings of the appropriate $w_{i_j}$s.  The ones to be checked should be obvious from Figure \ref{CombBigon}.  
$ $
$ $
	If the current pair neither intersects nor splits, we must trace further, advancing or decreasing the indices of the word segments in the same manner we used to get to the current pair of word segments(i.e. if we were looking at $W_{k-1}$ and $W_{l+1}$ for example, we must next look at $W_{k-2}$ and $W_{l+2}$.).  Continue checking in the given direction until an intersection is encountered, or no intersection is encountered but the segments split.  In the first case, we again have found a bigon.  In the second case, we must go back to our original pair of word segments $W_k$ and $W_l$ and check in a different direction, until all directions are exhausted.  If no direction results in finding another intersection, we conclude that our initial intersection $W_k$I$W_l$ cannot be a vertex of a bigon.  
	
	Starting with the pair P,C, we check all pairs of word segments for potential vertices of bigons until the list of pairs is exhausted or a vertex is found.  If no intersections are found then clearly our representative is minimal.  If an intersection is found, the procedure described above determines whether the vertex is part of a bigon or not.  If it is not, we continue checking the remaining pairs of word segments.  If the procedure finds a bigon, the method above gives us two sequences of word segments that were checked, starting with the two segments determining the first intersection and ending with the two word segments determining the second intersection.  There are four "orientations" associated with this pair of sequences, depending on if we moved in the forwards orientation for both $W_k$ and $W_l$ (this corresponds to the next pair of word segments being $W_{k+1},W_{l+1}$), forwards for the first and reverse for the second (i.e. the next pair of word segments are $W_{k+1},W_{l-1}$) and so on. These are denoted $(+,+)$, $(+,-)$,$(-,+)$, and $(-,-)$.  Note, however, that any pair of sequences with orientation $(-,-)$ starting with $W_k$ and $W_l$ can be written as a pair of sequences with orientation $(+,+)$ and ending with $W_k$ and $W_l$.  Similarly, a pair of sequences with orientation $(-,+)$ may be considered to be of the form $(+,-)$ simply by switching the order of the two sequences.  Thus, any pair of sequences determining a bigon essentially has orientation $(+,+)$ or $(+,-)$.  
	
	After a possible relabeling of the indices to fit with one of the essential orientations, we may abbreviate these two pairs of sequences with the notation $W_{k...i},W_{l...j}$ and call this pair of sequences a \textbf{combinatorial bigon}.  Each individual sequence is called a \textbf{combinatorial leg}.  The number $i-k mod n + 1 = j-l mod n + 1$ is called the \textbf{length} of the bigon.
	
	We summarize the preceding discussion with the following proposition:
	
	\start{prop}{bigons}{Bigons in a representative are in one-to-one correspondence with combinatorial bigons}\end{prop}
	
\subsection{Removing Bigons}	
\label{sec:rb}	
\begin{figure}[http]
\center{\includegraphics[width=3.0in]{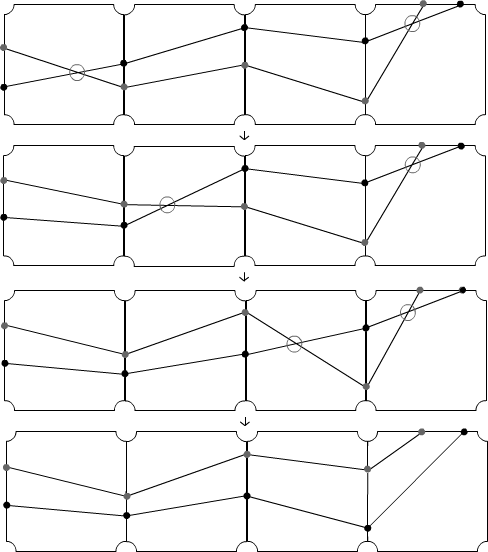}}
\caption{Geometrically, we are homotoping the bigon away.  Combinatorially, we are switching the locations of the grey and black points in P.}
\label{longbigon}
\end{figure}

Once we have a bigon, we would like to remove it as illustrated in Figure \ref{longbigon}.  The homotopy which removes the bigon of the representative induces a permutation on the elements in P.  It is easy to see what the permutations are if we break the homotopy into separate steps as in Figure \ref{longbigon}.  Conversely, certain permutations of points in P correspond to a homotopy in our actual representative.  Given a combinatorial bigon $W_{k...i},W_{l...j}$, we look at the pairs of word segments $W_k,W_l$, followed by $W_{k+1}$,$W_{l \pm 1}$, followed by $W_{k+2},W_{l \pm 2}$, ... ,$W_i,W_j$ (with the sign of $\pm$ depending on whether orientation is $(+,+)$ or $(+,-)$) and within each pair of word segments, choose the appropriate points to permute.  Which points to permute clearly depends on which orientation our combinatorial bigon has:
$ $

For $m \in \{1,2,...,L-1\}$ where $L$ is the length of the bigon.
\begin{enumerate}
\item	  $(+,+)$: switch $w_{{(k+m)}_2}$ with $w_{{(l+m)}_2}$ and $w_{{(k+m)}_1}$ with $w_{{(l+m)}_1}$
\item	  $(+,-)$: switch $w_{{(k+m)}_2}$ with $w_{{(l+m)}_1}$ and $w_{{(k+m)}_1}$ with $w_{{(l+m)}_2}$

\end{enumerate}

This does not remove self intersections in all cases.  There are certain conditions on the combinatorial bigon $W_{k...i},W_{l...j}$ that determine whether the above permutations remove intersections or not.  If no $W_{k+m}$ equals any $W_{l+m'}$, $m,m'$ $\in$ $\{1,2,...,L-1\}$ (i.e., the two combinatorial legs contain no word segments in common), then we have the following important theorem:

\start{theo}{defproper}{If the two combinatorial legs of a bigon $W_{k...i},W_{l...j}$ contain no word segments in common, then the permutation procedure removes the initial pair of intersections of the bigon and creates no additional intersections}\end{theo}

\begin{proof}

We first note that each pair of intermediate segments $W_{k + 1}W_{l /pm 1}$,...,$(W_{i - 1}W_{j /pm 1}$ becomes crossed, and then uncrossed as the permutations run their course, and in the end all of the intermediate word segments will be swapped(see Figure \ref{longbigon}). This clearly contributes no additional intersections.  Since each word segment is unique, only one point in each of of the terminal word segments $W_k$,$W_l$,$W_i$, and $W_j$ is permuted. After performing these permutations, one also sees that $W_k$ no longer intersects $W_l$, and likewise $W_i$ no longer intersects $W_j$.  

The only thing left to be checked is that no new intersections are formed during this process.  To see this, consider Figure \ref{defprop}
where $(a,b)$ and $(c,d)$ represent $W_k$ and $W_l$. Without loss of generality, we assume our bigon has orientation $(+,+)$ since the argument is nearly identical for the $(+,-)$ case.
   
\begin{figure}[http]
\center{\includegraphics[width=3.0in]{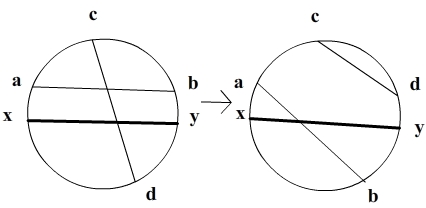}}
\caption{}
\label{defprop}
\end{figure}

Suppose a segment $(x,y) \neq (c,d)$  intersects $(a,b)$ after permuting the appropriate points, but not before.  Since $(x,y)$ did not intersect $(a,b)$ before permuting, both x and y must have been in a single arc with endpoints a and b. Since $(x,y)$ intersects $(a,b)$ afterwards, we may assume without loss of generality, we see that x and y must each have been in a separate arc determined by b and d before the permutations.  Thus $(x,y)$ must have intersected $(c,d)$ before permuting.  But at once we see that after permuting b with d, $(x,y)$ no longer intersects $(c,d)$.  Thus, for every intersection involving $(a,b)$ that we create, we remove an intersection with $(c,d)$.  We arrive at a symmetric statement for intersections involving $(c,d)$ and other word segments.  A nearly identical argument works using the segments $W_i$ and $W_j$, and so to avoid repetition we leave it to the curious.  Combining the facts we have collected, we arrive at the desired proposition.

\end{proof}

\begin{figure}[width=.8in][http]
\center{\includegraphics{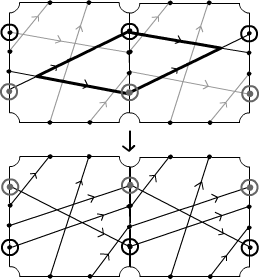}}
\caption{The surface is abAB, and the reduced word is bbAAAA.  Switching the indicated points does reduce the number of intersections.}
\label{improp2}
\end{figure}

If one or more word segments are shared, then the permutations just described may or may not reduce the number of intersections.  For instance, consider the situation depicted in Figure \ref{improp2}.  The indicated combinatorial legs share a word segment, and upon permuting the points as described above, we do not actually remove any intersections.  We call such combinatorial bigons \textit{non-removable}.  It will later be shown that non-removable combinatorial bigons correspond to improper bigons in the representative.

Before proceeding further, we present the following propositions which characterize the combinatorial symmetry that bigons with one or more equal segments in both combinatorial legs must have.

\start{prop}{shared} \begin{itemize}
\item (a) If the combinatorial bigon  $W_{k...i},W_{l...j}$ has one word segment shared by both combinatorial legs, it must be the first word segment of one leg, and the last word segment of the other.(i.e. $W_k=W_j$ or $W_i=W_l$)
$ $

\item (b)  If there are two word segments shared by both combinatorial legs, then they must either be the first two segments of one combinatorial leg and the last two of the other, or the first and last of both.
$ $

\item (c)  If there are $\geq$ 3 word segments shared by both combinatorial legs, then at least the first two word segments of one leg must be equal to the last two of the other leg.
\end{itemize}
\end{prop}

\begin{proof}
(a) and (b): Let one leg of a geometric bigon corresponding to our combinatorial bigon be denoted $L_1$ and the other $L_2$. Let $f$ $S \rightarrow M$ be the map of the curve.  We claim that the intersection $I$ of the preimages $f^{-1}(L_1)$ and $f^{-1}(L_2)$ can only be a proper subset of both $f^{-1}(L_1)$ and $f^{-1}(L_2)$; Indeed if $I = f^{-1}(L_1) = f^{-1}(L_2)$, then both legs of the bigon are identical, which cannot happen in our construction; if  $I = f^{-1}(L_1)$ $\neq$ $f^{-1}(L_2)$, then $L_1$ $\subset$ $L_2$ which implies both legs cannot have equal combinatorial length, a contradiction. Analogously for $I = f^{-1}(L_2)$ $\neq$ $f^{-1}(L_1)$.  Since the preimages must be connected, we have a situation similar to that depicted in Figure \ref{messprop1} , and only the ends of the preimages may intersect.  

\begin{figure}[http]
\center{\includegraphics[width=1.5in]{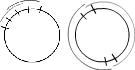}}
\caption{}
\label{messprop1}
\end{figure}

The proposition follows for cases (a) and (b).  For case (c), suppose for contradiction that no two shared word segments are adjacent in the combinatorial legs.  Choose one word segment shared by both legs. The preimage of the associated geometric leg, $f^{-1}(L_1)$ then, must wrap around the circle and intersect the preimages of the line segments determined by at least 3 word segments.  $f^{-1}(L_2)$ must wrap around the complementary part of the circle (otherwise we would see that at least two word segments are shared and adjacent) and intersect the same triple of preimages.  But this immediately gives rise to a contradiction.
\end{proof}

Now we show which types of combinatorial bigons always give rise to improper bigons in any representative, regardless of how the points are spaced on the fundamental polygon.  By Theorem \ref{hs} we can then simply skip all combinatorial bigons of these types when searching for removable bigons.  It happens that only a relatively small number of essential cases remain, and it will be shown that the permutations strictly reduce the number of intersections in these situations.  By Theorem \ref{defproper}, non-removable bigons must have at least one word segment shared by both combinatorial legs.  Furthermore, any such bigon must have orientation $(+,+)$, for otherwise, the segments shared by both are seen to inherit opposite orientations at the same time - clearly impossible. 

\start{theo}{conditions2}  Let $W_{k...i},W_{l...j}$ be a combinatorial bigon.  If one of the following conditions holds, then the corresponding geometric bigon in any segmented representative with the same P and C is improper:
\begin{enumerate}
\item At least two word segments are shared by both combinatorial legs

\item One segment be shared by both combinatorial legs, say $W_k=W_j$ and $w_{l_1}$ is between $w_{i_1}$ and $w_{k_1}$, and $w_{i_2}$ is between $w_{l_2}$ and $w_{k_2}$ in P.
\end{enumerate}

\end{theo}

\begin{proof}
First choose a segmented representative and a fundamental polygon based on on $P$ and $C$.  As in Proposition \ref{shared}, let one leg of the geometric bigon be denoted $L_1$ and the other $L_2$.  The preimages of the points where our representative crosses the edges of the polygon(i.e. the preimages of the points $w_{k_1},w_{k_2},w_{l_1},w_{l_2}$, etc.) partition the circle into N arcs, where N is the length of the reduced cyclic word used to construct our representative.  Suppose that there are at least two adjacent word segments in each combinatorial leg.  This is equivalent to saying that the intersection $I = f^{-1}(L_1) \cap f^{-1}(L_2)$ is not entirely contained in a single partitioning arc.  The only deformations of our segmented representative allowed are sliding the endpoints of segments along the edges of the fundamental polygon without changing P.  Any such homotopy of the line segment positions obviously leaves the combinatorial description of the bigon unchanged. Therefore we still have that $I$ is not contained within a single partitioning arc.  With our construction, the endpoints of the preimages of $L_1$ and $L_2$ must be contained in the interiors of the dividing arcs, for if they were not, then we would have two line segments emanating from the same point on an edge of the polygon, contrary to the endpoints of each line segment being distinct.  Since the endpoints of $f^{-1}(L_1)$ and $f^{-1}(L_2)$ must be contained in the interiors of dividing arcs, and since $I$ is not contained within a single dividing arc, $I$ cannot be empty.  Figure \ref{messprop1} shows this situation.  This covers all cases except when exactly two word segments are shared and that they occur at the ends of both combinatorial legs i.e. $W_k=W_j$ and $W_i=W_l$.  From this we see that our bigon really only has one vertex, and thus that our bigon must be improper.

For the 2nd condition above, our representative essentially looks like Figure \ref{essimprop2}.  The only way to remove the overlap is to translate one line segment over another, which clearly changes P.  Thus our combinatorial condition guarantees that the bigon will be improper.

\begin{figure}[http]
\center{\includegraphics[width=1.8in]{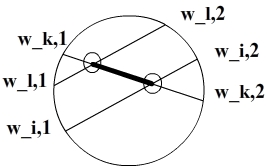}}
\caption{}
\label{essimprop2}
\end{figure}

\end{proof}

\start{theo}{conditions1} Let $W_{k...i},W_{l...j}$ be a combinatorial bigon with $W_k=W_j$, and no other segments shared. Then the bigon is removable if the 2nd condition from Theorem \ref{conditions2} is not satisfied.
\end{theo}

\begin{proof}
As in Theorem \ref{defproper} we only need to consider the 6 points in the terminal word segments, $W_k$,$W_l$, and $W_i$, since the intermediate segments are swapped. We consider all possible relative positions of those 6 points within P that simultaneously realize (1)$(w_{k_1},w_{k_2})I(w_{l_1},w_{l_2})$, (2) $(w_{k_1},w_{k_2})I(w_{i_1},w_{i_2})$, (3) $w_{k_2}$ and $w_{l_2}$ on the same edge, and (4) $w_{k_1}$ and $w_{i_1}$ on the same edge. These four conditions are necessary in order for our pair of sequences to actually be a combinatorial bigon. There are 12 configurations of points satisfying these conditions, and only 6 that need be considered once symmetry is taken into account. To see this, consider an oriented circle representing the cyclic ordering of P. We imagine placing the initial line segment, say $(w_{k_1},w_{k_2})$ which then divides the circle into two arcs. Next, there are two choices for the placement of the line segment $(w_{l_1},w_{l_2})$, corresponding to which arc of the circle we wish to place the point $w_{l_1}$. For a given placement of $(w_{l_1},w_{l_2})$, the circle is divided into four arcs. The point $w_{i_1}$ may be placed on any arc except the one between the points $w_{k_2}$ and $w_{l_2}$ (that would imply $w_{i_1}$ is not on the same edge as $w_{k_1}$ or that $w_{k_2}$ and $w_{l_2}$ are on the same edge, contrary to conditions (3) or (4). Once $w_{i_1}$ is placed, there are two choices for $w_{i_2}$ that give $(w_{k_1},w_{k_2})I(w_{i_1},w_{i_2})$. Counting all combinations, we get $2*3*2=12$ total, but we can eliminate half of the permutations, since one ordering and its reverse ordering are equivalent for the purpose of determining the effect of the permutations. 

Below, the 6 essential orderings and total number of intersections they determine are given, along with their ordering and intersections after the permutations are performed:

 1) $w_{k_1}$-$w_{i_1}$-$w_{l_2}$-$w_{k_2}$-$w_{i_2}$-$w_{l_1}$ 3 $ $ $w_{i_1}$-$w_{k_1}$-$w_{k_2}$-$w_{l_2}$-$w_{i_2}$-$w_{l_1}$ 1

$ $
2) $w_{k_1}$-$w_{i_1}$-$w_{l_2}$-$w_{k_2}$-$w_{l_1}$-$w_{i_2}$ 2 $ $ $w_{i_1}$-$w_{k_1}$-$w_{k_2}$-$w_{l_2}$-$w_{l_1}$-$w_{i_2}$ 0

$ $
3) $w_{k_1}$-$w_{l_2}$-$w_{i_2}$-$w_{k_2}$-$w_{l_1}$-$w_{i_1}$ 3 $ $ $w_{i_1}$-$w_{k_2}$-$w_{i_2}$-$w_{l_2}$-$w_{l_1}$-$w_{k_1}$ 1

$ $
4) $w_{k_1}$-$w_{i_2}$-$w_{l_2}$-$w_{k_2}$-$w_{l_1}$-$w_{i_1}$ 2 $ $ $w_{i_1}$-$w_{i_2}$-$w_{k_2}$-$w_{l_2}$-$w_{l_1}$-$w_{k_1}$ 0

$ $
5) $w_{k_1}$-$w_{i_2}$-$w_{l_2}$-$w_{k_2}$-$w_{i_1}$-$w_{l_1}$ 3 $ $ $w_{i_1}$-$w_{i_2}$-$w_{k_2}$-$w_{l_2}$-$w_{k_1}$-$w_{l_1}$ 1

$ $
6) $w_{k_1}$-$w_{l_2}$-$w_{i_2}$-$w_{k_2}$-$w_{i_1}$-$w_{l_1}$ 2 $ $ $w_{i_1}$-$w_{k_2}$-$w_{i_2}$-$w_{l_2}$-$w_{k_1}$-$w_{l_1}$ 2

Case 6 satisfies the second condition in Theorem \ref{conditions2} and so determines an improper bigon, so it is not surprising that the permutations do not reduce intersections.  

The only thing that remains to be proved is that for each of the 5 "good" cases, no additional intersections are produced.  To make things easier to keep track of, let us relabel the word segments as follows: $W_k=(a,b)$,$W_l=(c,d)$, and $W_i=(e,f)$. 

We now prove the theorem for case 1 above by following a line of reasoning similar to that in Theorem \ref{defproper}.  Careful consideration of Figure \ref{c1pf} yields the proof, but we give some of the details here.  Suppose a word segment $(x,y) \neq (a,b),(c,d)$ or $(e,f)$, does not intersect $(a,b)$ initially but does after permuting.  Then $(x,y)$ must have intersected one of $(e,f)$ or $(c,d)$ beforehand.  To see this, note that one of the points x or y must lie in the "top" arc between e and d, while the other point must have been either between a and e, or d and b.  However, after the permutations, it cannot intersect either of them, as seen from the figure and the argument in the previous sentence.  Now suppose $(x,y)$ intersects $(c,d)$ after permuting but not before.  Since the point c is unchanged by a permutation, we see that $(x,y)$ must have intersected both $(a,b)$ and $(e,f)$ before the permutations.  But clearly, after the permutations $(x,y)$ could not possibly intersect $(a,b)$, since then $(x,y)$ would have intersected $(c,d)$ to begin with.  Finally, suppose $(x,y)$ intersects $(e,f)$ after the permutations but not before.  The point f is fixed by the permutations, and we see that $(x,y)$ must have intersected both $(a,b)$ and $(c,d)$ before the permutations.  By the exact same reasoning as the previous case, $(x,y)$ can no longer intersect $(a,b)$ after the permutations.  In every situation, if an intersection is introduced, there is a corresponding intersection that is removed.  Thus there is no net gain of intersections between $(x,y)$ and the word segments $(a,b),(c,d)$ and $(e,f)$.  Since the word segment $(x,y)$ was arbitrary, we proved the theorem for case 1.  The proofs for the other 4 "good" cases are nearly identical, so we leave them to the curious.

\begin{figure}[http]
\center{\includegraphics[width=3.0in]{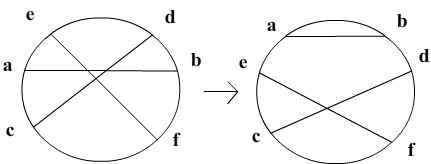}}
\caption{}
\label{c1pf}
\end{figure}

\end{proof}

\section{Summary}

Once given a surface and a reduced word, we construct P and C in whatever appropriate manner we prefer.  Then we examine P and C to determine if there are any combinatorial bigons.  If there are none, then by Proposition \ref{bigons} and Theorem \ref{hs} our current representative is minimal, and we are done.  If we find a bigon, but it satisfies one of the conditions Theorem \ref{conditions2}, then we skip over it and search for more bigons.  If every bigon is of one of the types described in Theorem \ref{conditions2} then again we are done, by the earlier propositions and Theorem \ref{hs}.  Finally, if we find a bigon not of the types in Theorem \ref{conditions2}, then by Theorem \ref{conditions1} we know we may apply the permutations and strictly reduce the total number of intersections encoded in P and C.  This process must terminate, and we eventually end up with a modified P and C which encodes a minimal segmented representative.

\section{Example}

\begin{figure}[http]
\center{\includegraphics[width=3.0in]{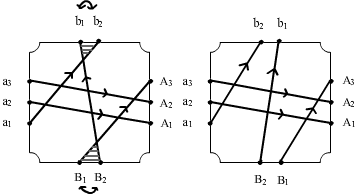}}
\caption{}
\label{alg3}
\end{figure}

We will now show how this algorithm works for the surface word \textbf{abAB} and the reduced cyclic word \textbf{bbAAA}, the same surface and curve as in Figure \ref{Alg1}.  Recall that in Figure \ref{Alg1}, we have chosen P = $a_1, a_2, a_3, b_1, b_2, A_3, A_2, A_1, B_2, B_1$ and C = $(B_1,A_3) (a_3,A_2) (a_2,A_1) (a_1,b_2) (B_2,b_1)$ as our initial representative.  We now check through pairs of word segments until we find a pair that intersect.  Upon inspection, we see that $(a_3,A_2)I(B_2,b_1)$.  The next step is to see if this vertex can possibly be part of a bigon.  There can be no bigon starting at this vertex with a $(+,+)$ orientation, since the segments spit to the A and b sides.  Likewise for the other 3 orientations.  We conclude that this particular intersection cannot be a vertex of a bigon and move on until we find another pair of segments that intersect.  Suppose the next pair of segments we find to intersect are $(B_2,b_1)$ and $(a_1,b_2)$.  For this particular intersection, we see that it may be a vertex of a bigon with a $(+,+)$ orientation, since the points $b_1$ and $b_2$ are on the same edge.  The next pair of segments we compare are $(B_1,A_3)$ and $(B_2,b_1)$, which also intersect.  Thus we have a combinatorial bigon \{$(B_2,b_1)(B_1,A_3);(a_1,b_2),(B_2,b_1)$\}.  Since $(B_2,b_1)$ is shared by both sequences, we must use Theorem \ref{conditions1} to check if we can remove this bigon to reduce the number of self intersections.  We see that this particular bigon is of type (5) in reverse, and so we may proceed.  We need to switch $b_1$ with $b_2$ and $B_1$ with $B_2$, as shown in Figure \ref{alg3}.

\begin{figure}[http]
\center{\includegraphics[width=3.0in]{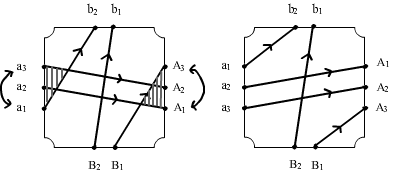}}
\caption{}
\label{alg4}
\end{figure}

Now, we have P = $a_1, a_2, a_3, b_2, b_1, A_3, A_2, A_1, B_1, B_2$, while C remains the same.  We must again compare pairs of segments until we find an intersection.  We see that $(a_1,b_2)I(a_2,A_1)$, so we check to see if this intersection can be a vertex of a bigon.  A $(-,-)$ orientation is ruled out, since $b_2$ and $A_1$ are on different edges, but we see that we may check for a $(+,+)$ orientation.  Doing so results in the combinatorial bigon \{$(a_1,b_2),(a_2,A_1),(a_3,A_2);(a_2,A_1),(a_3,A_2),(B_1,A_3)($\}.  This bigon has two segments shared by both sequences, so it is guaranteed to be improper by Theorem \ref{conditions2}.  By Theorem \ref{hs} there must be a proper bigon if our representative does not have minimal self-intersection, so we continue to look for a different bigon to potentially remove.  Suppose we next find that $(a_1,b_2)I(a_3,A_2)$.  We find that it is the vertex of the bigon \{$(a_1,b_2),(a_2,A_1);(a_3,A_2),(B_1,A_3)$\}.  We switch $a_1$ with $a_3$, and $A_1$ with $A_3$, as shown in Figure \ref{alg4}.  Finally, we once again check all pairs of segments using the permuted P, and determine that there are no more proper bigons, and thus by Theorem \ref{hs} our representative has the minimal number of self intersections possible.  The final output of the algorithm is P = $a_3, a_2, a_1, b_2, b_1, A_1, A_2, A_3, B_1, B_2$ and C = $(B_1,A_3) (a_3,A_2) (a_2,A_1) (a_1,b_2) (B_2,b_1)$, which determines our representative.

\section{Remarks}

This algorithm provides an alternative way to combinatorially count the minimal self intersections of a free homotopy class on a surface.  Its efficiency is polynomial in the length of the curve word, where the "units" used are comparisons and permutations.  The concepts presented here will also work for finding a minimally intersecting configuration of two or more curves and only need to be slightly modified.  One surprise that was encountered during this work was that the combinatorial description of a bigon alone is not always sufficient to determine if the corresponding geometric bigon is improper or not.  There are combinatorial bigons which may or may not correspond to improper bigons depending on the actual positioning of the points along the boundary of the fundamental polygon. This is why we use the words "removable" and "non-removable" when dealing with combinatorial bigons.

$ $

\bibliographystyle{amsalpha}

\begin{thebibliography}{77} 
%
\bibitem{BS} J.\ Birman and C.\ Series, An algorithm for simple curves on surfaces, {\em J. London Math. Soc.} (2), {\bf 29} (1984), 331-342.
\bibitem{chas1} M.\ Chas, Combinatorial Lie Bialgebras of curves on surfaces;\href{http://arxiv.org/abs/math/0105178v2}{\tt arXiv: 0105178v2 [math.GT]}
\bibitem{chas2} M.\ Chas, Minimal intersection of curves on surfaces; arXiv:0706.2439v2 [math. GT]%
\bibitem{cl} M.\ Cohen and M.\ Lustig, Paths of geodesics and geometric intersection numbers I,  {\em Combinatorial Group Theory and Topology}, Alta, Utah, 1984,  Ann. of Math. Studies  {\bf 111}, Princeton Univ. Press, Princeton,  (1987), 479-500.
\bibitem{hassscott} J. Hass, P. Scott, Intersections of curves on surfaces, Israel J. Math. Vol 51 (1985), pgs 90-120
\bibitem{still} J. Stillwell, Classical Topology and Combinatorial Group Theory, 2nd Ed., Springer-Verlag, 1993
%
%
%


%
\end{thebibliography}

\end{document}